\documentclass[11pt, a4paper]{article}
\usepackage{latexsym}
\usepackage{indentfirst}
\usepackage{graphicx}
\usepackage{amsmath}
\usepackage{amssymb}
\begin{document}
\title{New results concerning the stability of equilibria of a
delay differential equation modeling leukemia}
\author{Anca-Veronica Ion\\
"Gh. Mihoc-C. Iacob" Institute of Mathematical Statistics\\
and Applied Mathematics of the Romanian Academy\\ Bucharest, Romania
}

\date{}
\maketitle

\begin{abstract}
The paper is devoted to the study of stability of equilibria of a
delay differential equation that models leukemia. The equation was
previously studied in  [5] and  [6], where the emphasis is put on
the numerical study of periodic solutions. Some stability results
for the equilibria are also presented in these works, but they are
incomplete and contain some errors. Our work aims to complete and to
bring corrections to  those results. Both Lyapunov first order
approximation method and second Lyapunov method
are used.\\
\textbf{Acknowledgement.} Work supported by Grant 11/05.06.2009
within the framework of the Russian Foundation for Basic Research -
Romanian Academy collaboration.
\end{abstract}

\section{Introduction}
The study of the mathematical model of periodic chronic myelogenous
leukemia considered in \cite{PM-M}, \cite{PM-B-M} may be reduced to
that of the equation
\begin{equation}\label{eq}
\dot{x}(t)=-\left[\frac{\beta_0}{1+x(t)^n}+\delta\right]x(t)+k\frac{\beta_0x(t-r)}{1+x(t-r)^n},
\end{equation}
where $\beta_0,\,n,\,\delta,\,k,\,r$ are positive parameters and
$k=2e^{-\gamma r},$ with $\gamma>0$.  We do not insist here on the
significance of the function $x(.)$ or in that of the parameters,
since these are extensively presented in \cite{PM-M}, \cite{PM-B-M}.
We only remind that the unknown function, $x(\cdot),$ should be
nonnegative, being a non-dimensional density of cells.

The two studies cited above are mainly devoted to the numerical
investigation of the delay differential equation above. The
 stability of equilibria study in \cite{PM-M}, \cite{PM-B-M}, that is reduced
to a few lines, is, unhappily, incomplete and contains  some errors.

Our work aims to correct the errors in the stability conditions
presented in \cite{PM-M}, \cite{PM-B-M} and to present some aspects
concerning the dynamics generated by equation (\ref{eq}), aspects
that were not pointed out there.

In Subsection 1.1  we prove that the Cauchy problem associated to
our equation has an unique defined on $[-r,\infty)$ bounded
solution. In Subsection 1.2, following \cite{PM-M}, \cite{PM-B-M},
the two equilibrium solutions, as well as the linearized equation
and the characteristic equation for each of these, are presented.
Section 2 deals with the stability study of the equilibrium points.
We use the results of \cite{Ha} and, in a first stage, we obtain
results valid for both equilibria. In the subsequent subsections we
analyze the stability of the two points individually. It is
important to perform this separate study since the conclusions are
very specific to each equilibrium point. In Section 3 we comment the
stability results in \cite{PM-M}, \cite{PM-B-M}, pointing out the
errors therein.

\subsection{Existence and uniqueness of solution}
 We  make the notation
 $\mathcal{B}=C([-r,0],\mathbb{R})$ (the space of continuous, real valued functions
 defined on $[-r,0],$ with the supremum norm, denoted by $|x|_0$).
Given a function $x:[-r,T)\mapsto \mathbb{R},\,T>0$ and a $0\leq t <
T,$ we define the function $x_t\in \mathcal{B}$ by $x_t(s)=x(t+s).$

Equation (\ref{eq}) may be written as \begin{equation}\label{abseq}
\dot{x}=h(x_t),
\end{equation}
where $h:\mathcal{B}\mapsto\mathbb{R},$ and  we impose to this
equation the initial condition
\begin{equation}\label{ci}x_0=\phi\in \mathcal{B}.
\end{equation}

Remark that if the initial condition is a positive function, then
$x(t)$ can not become strictly negative. Indeed, let $t_1$ be the
first moment when $x(t_1)=0,$ (that is $x(t)>0$ for $t<t_1$). Then
$\dot{x}(t_1)=k\frac{\beta_0x(t_1-r)}{1+x(t_1-r)^n}>0,$ hence
$x(t)>0$ for $t\geq t_1$ in a neighborhood of $t_1$.

We study the existence, uniqueness and domain of existence of
solutions.

The function $h$ is globally Lipschitz. Indeed, by denoting
\\$\beta(x)=\beta_0/(1+x^n),$ we have
\[|\frac{d}{dx}(\beta(x)x)|<\beta_0(n+1),
\]
and thus, for any $\varphi_1,\,\varphi_2\in \mathcal{B}$
\[|h(\varphi_1)-h(\varphi_2)|\leq
(\beta_0(n+1)+\delta)|\varphi_1(0)-\varphi_2(0)|+k\beta_0(n+1)|\varphi_1(-r)-\varphi_2(-r)|\leq\]
\[\leq [(k+1)\beta_0(n+1)+\delta]|\varphi_1-\varphi_2|_0.
\]
From here the continuity of $h$ follows also.

Theorem 2.3 of \cite{HaL} implies that problem (\ref{eq}),
(\ref{ci}) has an unique solution defined on an interval $[0,\,T).$

If we take $\varphi_2=0,$ by using the fact that $h(0)=0,$ we obtain
that
\[|h(\varphi_1)|\leq [(k+1)\beta_0(n+1)+\delta]|\varphi_1|_0
\]
thus the function $h$ is also completely continuous.

We prove that the solution is bounded. For this, we multiply
equation (\ref{eq}) by $x(t)$ and we get
\[\dot{x}(t)x(t)=-\beta_0\frac{x^2(t)}{1+x^n(t)}-\delta
x^2(t)+k\beta_0\frac{x(t-r)x(t)}{1+x^n(t-r)}\leq
\]
\[\leq -\delta
x^2(t)+\frac{\varepsilon}{2}k\beta_0\frac{x^2(t)}{1+x^n(t-r)}+\frac{1}{2\varepsilon}k\beta_0\frac{x^2(t-r)}{1+x^n(t-r)}
\]
\[\leq -\delta
x^2(t)+\frac{\varepsilon}{2}k\beta_0
x^2(t)+\frac{1}{2\varepsilon}k\beta_0,
\]
and
\[\frac{d(x^2(t))}{dt} + (2\delta-\varepsilon
k\beta_0)x^2(t)\leq \frac{k\beta_0}{\varepsilon}.
\]
We chose an $\varepsilon>0$ such that $\eta:=2\delta-\varepsilon
k\beta_0>0$ and we obtain by integration \[x^2(t)\leq
\phi^2(0)e^{-\eta t}+\frac{k\beta_0}{\varepsilon \eta },\] hence the
solution of problem (\ref{eq}), (\ref{ci}) is bounded. Theorem 3.2
of \cite{HaL} implies that the solution is defined on the whole
positive real time semiaxis.

Hence, for any $\phi \in \mathcal{B},$ problem (\ref{eq}),
(\ref{ci}) has an unique defined on $\mathbb{R}^+$ bounded solution.
We can thus associate to this problem the semigroup of operators on
$\mathcal{B},$ $\{T(t)\}_{t\geq 0},$ $T(t)(\phi)=x_t(\phi)$, where
$x(t,\phi)$ is the solution of eq. (\ref{eq}) with initial condition
$x_0=\phi$.

\subsection{Equilibrium solutions}

In this subsection we, inevitably, follow  \cite{PM-M}.

The equilibrium points of the problem are $$x_1=0,\,\,
x_2=(\frac{\beta_0}{\delta}(k-1)-1)^{1/n}.$$ The second one is
acceptable from the biological point of view if and only if it is
strictly positive that is, if and only if
\begin{equation}\label{cond} \frac{\beta_0}{\delta}(k-1)-1
> 0.\end{equation}

In terms of $r,$ by using $k=2e^{-\gamma r},$ the above inequality
may be written as
\begin{equation}\label{rmax}  r<r_{max}:=-\frac{1}{\gamma}\ln\frac{1}{2}\left(1+\frac{\delta}{\beta_0}\right),
\end{equation}
and since the delay $r$ is positive, the condition
$\delta/\beta_0<1$ follows.

The biological interpretation of function $\beta$ \cite{PM-M} shows
that the condition $\beta(x_2)=\delta/(k-1)>0$ should be fulfilled.
This is equivalent to $k>1.$

The linearized equation around one of the equilibrium points is
\begin{equation}\label{lineq}
\dot{z}(t)=-[B+\delta]z(t)+kBz(t-r),\end{equation} with
$B=\beta'(x^*)x^*+\beta(x^*),\,x^*=x_1$ or $x^*=x_2.$

The eigenvalues of the infinitesimal generator of the semigroup of
operators generated by equation (\ref{lineq}) are the solutions of
the characteristic equation
\begin{equation}\label{chareq}
\lambda+\delta+B=kBe^{-\lambda r}.
\end{equation}

\section{Stability of equilibrium points}

In order to investigate the stability of the equilibrium solutions,
we first try to establish the conditions in which all the
eigenvalues have strictly negative real part, in order to use the
linear approximation Lyapunov method.

We rely on the work \cite{Ha} that exhaustively solves the problem
of finding necessary and sufficient conditions on the parameters
such that the equation $\lambda=a_1+a_2e^{-\lambda}$ has only
solutions with strictly negative real part.

We denote $\delta+B=p,\,kB=q,$ hence (\ref{chareq}) becomes
\begin{equation}\label{chareq1}
\lambda+p=qe^{-\lambda r}.\end{equation} By taking
$\lambda=\mu+i\omega,$ and by equating the real, resp. the imaginary
parts in our equation, we obtain
\begin{equation}\label{re,im}
\mu+p=q e^{-\mu r}\cos(\omega r),
\end{equation}
\[\omega=-q e^{-\mu r}\sin(\omega r).
\]
It is useful to consider the case $\mu=0$ in the above equations,
\begin{equation}\label{muegal0}
p=q \cos(\omega r),
\end{equation}
\[\omega=-q \sin(\omega r).
\]

The results in \cite{Ha} imply the following

\textbf{Proposition} \textit{All solutions $\lambda$ of eq.
(\ref{chareq1}) satisfy  $Re \lambda<0,$ if and only if\\
\textbf{a)} $q<0,\,0<-pr<1$ and
$-p<-q<(\omega_0^2+p^2)^{1/2},$\\or\\ \textbf{b)}
$q<0,\,p>0$ and $-q<(\omega_0^2+p^2)^{1/2},$\\or\\
\textbf{c)} $q>0,\,p>0,$ and $q<p,$\\
where $\omega_0$ is the solution in  $(0,\pi/r)$ of the equation}
\begin{equation}\label{omega-0}\omega\cot (\omega r)=-p.
\end{equation}

\textbf{Remark.} If we divide the first equality in (\ref{muegal0})
to the second one, we obtain (\ref{omega-0}). Hence relations
(\ref{muegal0}) are equivalent to the set of relations
(\ref{omega-0}) and $\omega^2+p^2=q^2.$
\bigskip

In order to express $\omega_0$ in a more direct form, we consider
the function $T:[0,\pi)\mapsto (-\infty,\,1],$ given by
\begin{equation}\label{T}
T(y)=\left\{%
\begin{array}{ll}
    y\cot(y), &  y\in (0,\pi); \\
    1, & y=0. \\
\end{array}%
\right.
\end{equation}
The function is a bijection and we can equivalently define
$\omega_0$, the solution of (\ref{omega-0}), as
\begin{equation}\label{om0}\omega_0=\frac{1}{r}T^{-1}(-pr).\end{equation}

We express the conditions in Proposition in terms of $r.$  We first
remark that
\[\omega_0^2+p^2=\frac{p^{2}}{\cot^2(\omega_0 r)}+p^2=\frac{p^2}{\cos^2(\omega_0
r)}.
\]
\\
The two inequalities in \textbf{ \textit{a)}} of the above
Proposition may be written as
\[|p|<|q|<\frac{|p|}{|\cos(\omega_0 r)|}.
\]

Since $p<0,$  the solution $\omega_0$ of equation (\ref{omega-0}) is
such that $\omega_0 r\in (0,\pi/2).$ Hence the above inequality is
equivalent to
\[0<\frac{p}{q}<1,\,\cos(\omega_0 r) <\frac{p}{q},
\]
and the second one is equivalent to $\arccos(\frac{p}{q})<\omega_0
r<\pi/2.$ To conclude, case \textbf{a)} is described by the
inequalities
\begin{equation}
q<p<0,\frac{\arccos(\frac{p}{q})}{\omega_0}<r< \frac{1}{|p|}.
\end{equation}

In case \textbf{b)} $q<0,\,p>0,$ and we must have
\[-q<\frac{p}{|\cos (\omega_0 r)|}.
\]
In this case, $\omega_0\cot \omega_0 r=-p<0,$ and since $\omega_0
r\in (0,\pi), $ we must have $\omega_0 r\in(\pi/2,\pi).$  The above
inequality is equivalent to
\[|\cos (\omega_0 r)|<\frac{p}{|q|}
\]
and this one is satisfied if  \[p/|q|>1\,\, \mathrm{or}\,\,
\{p/|q|\leq 1\,\,\mathrm{and}\, (-\cos (\omega_0
r)<\frac{p}{-q})\}.\] The last condition is equivalent to \[-1\leq
\frac{p}{q}< 0\,\mathrm{and}\, \cos (\omega_0
r)>\frac{p}{q}\,\Leftrightarrow \frac{\pi}{2}<\omega_0
r<\arccos(\frac{p}{q})\,\Leftrightarrow
\]
\[\Leftrightarrow\frac{\pi}{2\omega_0}<r<\frac{\arccos(\frac{p}{q})}{\omega_0}.
\]

\textbf{Remark.} The case $q<0,\,p=0$ corresponds to $\omega_0
r=\pi/2,$ and the eigenvalues lie to the left of the vertical axis
if and only if $-q r<\pi/2.$

We can now translate the discussion above to our concrete problem.\\
\textbf{I.} If $B<0,$ then two situations may occur.

\textbf{A.} $\delta+B<0.$ In this situation, $Re\lambda <0$ for all
eigenvalues $\lambda$ if and only if $|\delta+B|<|kB|$ and
\begin{equation}\label{A}
\frac{\arccos{((\delta+B)/kB)}}{\omega_0}<r<\frac{1}{|\delta+B|},
\end{equation}
where
\[\omega_0=\frac{1}{r}T^{-1}(-(\delta+B)r),
\]
with $T$ given by (\ref{T}).

If the studied equilibrium point is $x_2,$ the condition $r\leq
r_{max}$ must be also fulfilled.

\textbf{B.} $\delta+B>0.$ In this situation, $Re\lambda <0$ for all
eigenvalues $\lambda$ if and only if
\begin{equation}\label{B}
\delta+B>|kB|\,\,\mathrm{or}\,\left\{\delta+B\leq |kB|\,\,
\mathrm{and}\,\,r<\frac{\arccos{((\delta+B)/kB)}}{\omega_0}\right\}
\end{equation}
where, again \[\omega_0=\frac{1}{r}T^{-1}(-(\delta+B)r),
\]
with $T$ given by (\ref{T}).

\textbf{II.} If $B>0,$ then we can only have $\delta+B>0,$ and in
this situation $Re\lambda<0$ for all eigenvalues $\lambda$ if and
only if
\[kB<\delta+B.
\]

Even if the above discussion seems comprehensive, it is still useful
to consider the two equilibrium points separately and to discuss
their stability.
\subsection{Stability properties of $x_1$}
In this case, $B=\beta_0>0,$ hence the necessary and sufficient
condition for the negativity of the real part of all eigenvalues is
\[k\beta_0<\delta +\beta_0\Leftrightarrow\,\frac{\beta_0}{\delta}(k-1)<1.
\]
Since the condition (\ref{cond}) for the existence of the second
equilibrium point, $x_2,$ is the reverse of the above inequality, it
follows that $x_1$ is stable as long as it is the single equilibrium
point. When the second equilibrium point occurs, $x_1$ becomes
unstable.

We inspect the eigenvalues at $\frac{\beta_0}{\delta}(k-1)=1.$
Equation (\ref{chareq}) in this case is
\[\lambda+\delta+\beta_0=k\beta_0e^{-\lambda r},
\]
and, since $k\beta_0=\delta +\beta_0,$ admits the solution
$\lambda=0.$ Hence the change of stability occurs by traversing the
eigenvalue $\lambda=0.$

\subsubsection{Stability of $x_1$ when $\frac{\beta_0}{\delta}(k-1)=1$}

In this case, the "first order approximation" theorem is of no use,
since $0$ is the eigenvalue with greatest real part. We use a
Lyapunov function in order to prove stability of the zero solution.

However, since for our problem $x(t)\geq 0,$ the concept of
stability should be interpreted in the following way:\\
\textit{for every $\varepsilon >0$ there is a $\delta >0$ such that
if $\phi(s)\geq 0,\, s\in[-r,0],$ and $|\phi|_0<\delta$, then $0\leq
x(t,\phi)<\varepsilon$ for any $t>0,$ where, as above, $x(t,\phi)$
is the solution of (\ref{eq}) with condition (\ref{ci}).}

If $V:\mathcal{B}\mapsto\mathbb{R}$ is continuous, the derivative
along the solution $x(\cdot,\phi)$ of the Cauchy problem
(\ref{abseq}), (\ref{ci}) is defined as \cite{HaL}
\[\dot{V}(\phi)=\limsup_{h\rightarrow
0^{+}}\frac{1}{h}[V(x_{h}(\phi))-V(\phi)].
\]
\textbf{Definition} \cite{HaL}. \emph{V is a Lyapunov function on
$G\subset \mathcal{B}$} \textit{if} $V$ \textit{is continuous on
}$\overline{G}$ \textit{and} $\dot{V}\leq 0$ \textit{on G.}

 \textbf{Theorem} \cite{B}. \emph{If} $V:\mathcal{B}\mapsto\mathbb{R}$ \emph{is
 a Lyapunov function and there
 exist a continuous increasing function}
 $a:[0,\,\infty)\mapsto[0,\,\infty),$ \emph{with} $a(0)=0$ \emph{and}
\[a(|\phi(0)|)\leq V(\phi),
\]
\emph{then the solution }$x=0$ \emph{of equation (\ref{abseq}) is
stable and every solution is bounded.}

We construct below a Lyapunov function for our problem, for the
considered parameter values.

Let us consider the function $G(u)=\int_{0}^{u}2s/(1+s^n)ds.$  We
define
\[V(\phi)=G(\phi(0))+k\beta_0\int_{-r}^{0}\frac{\phi^{2}(s)}
{(1+\phi^{n}(s))^{2}}ds.
\]
We have
\[\dot{V}(\phi)= \frac{2\phi(0)}{1+\phi^{n}(0)}\dot{x}(0,\,\phi)+
k\beta_{0}\left[\frac{\phi^{2}(0)}
{(1+\phi^{n}(0))^{2}}-\frac{\phi^{2}(-r)}
{(1+\phi^{n}(-r))^{2}}\right],
\]
and by using the equality
\[\dot{x}(0,\phi)=-[\beta(\phi(0))+\delta]\phi(0)+k\beta(\phi(-r))\phi(-r),
\]
we obtain
\[\dot{V}(\phi)=-2\beta_0\frac{\phi^2(0)}{(1+\phi^{n}(0))^{2}}-
2\delta\frac{\phi^2(0)}{1+\phi^{n}(0)}+
2k\beta_0\frac{\phi(0)\phi(-r)}{(1+\phi^{n}(0))(1+\phi^{n}(-r))}+
\]
\[+k\beta_{0}\left[\frac{\phi^{2}(0)}
{(1+\phi^{n}(0))^{2}}-\frac{\phi^{2}(-r)}
{(1+\phi^{n}(-r))^{2}}\right],
\]
from where, with the inequality
\[\frac{2\phi(0)\phi(-r)}{(1+\phi^{n}(0))(1+\phi^{n}(-r))}\leq
\frac{\phi^2(0)}{(1+\phi^{n}(0))^{2}}+\frac{\phi^2(-r)}{(1+\phi^{n}(-r))^{2}}
\]
we obtain
\[\dot{V}(\phi)\leq 2(-\beta_0-\delta+k\beta_0)\frac{\phi^2(0)}{(1+\phi^{n}(0))^{2}}=0,
\]
since $k\beta_0=\delta +\beta_0.$

The hypotheses of Theorem 1 are satisfied with $a(u)=G(u),$ and it
follows that $x_1=0$ is stable in the case of this subsection.

\subsection{Stability properties of $x_2$}

In this case,
\begin{equation}\label{B1}
B=\beta_0[n-(n-1)A]/A^2\end{equation} where $A=\beta_0(k-1)/\delta.$

As pointed out in Subsection 1.2, in this case the condition
(\ref{rmax}) must be fulfilled.

We refine the discussion concerning the cases of stability given at
the beginning of Section 2, for this concrete $B$.

\textbf{I.A.} The condition $B<0$ and the definition of $B$ imply
$n-(n-1)A<0.$ This implies $n>1$ and
\begin{equation}\label{Bneg}\frac{\beta_0}{\delta}(k-1)>\frac{n}{n-1}.\end{equation}

The condition $B+\delta<0$ leads to $n-(n-k)A<0,$ that implies $n>k$
and
\begin{equation}\label{Bplusdneg}\frac{\beta_0}{\delta}(k-1)>\frac{n}{n-k}.\end{equation}
Obviously, the second inequality implies the first one.

The sufficient condition of local stability is condition (\ref{A}),
with $B$ given by (\ref{B1}). We remark that the condition
$|\delta+B|<|kB|$ is satisfied since it is equivalent to
$\delta+B>kB$ and this one is equivalent to
$\frac{\beta_0}{\delta}(k-1)>1,$ (the condition of positivity of
$x_2$).

We have to study the behavior of the solutions at the extremities of
the interval of stability.

\textbf{a)} We consider the case
\begin{equation}\label{r-H}
r=\frac{\arccos((\delta+B)/(kB))}{\omega_0}.
\end{equation}
This relation, together with
 \begin{equation}\label{om01}\omega_0\cot(\omega_0
r)=-(\delta+B),\end{equation} and $\omega_0>0$ (from the definition
of $\omega_0$), imply
 \begin{equation}\label{om02}\omega_0=\sqrt{(kB)^2-(\delta+B)^2}\end{equation} and that the pair
$\mu^*=0,\,\omega^*=\omega_0$ represents a solution of
(\ref{re,im}).

For later use we remark that, for $B<0,$ the relations \eqref{om01}
 and \eqref{om02} (where, by the definition of $\omega_0,$  $\omega_0 r\in (0,\pi)$) together, imply
relation \eqref{r-H} and again that the pair
$\mu^*=0,\,\omega^*=\omega_0$ is a solution of (\ref{re,im}).

We assume that we vary one of the parameters, that we denote here by
$\alpha$, such that for a value $\alpha^*$ the equality (\ref{r-H})
is satisfied, and keep all other parameters fixed. We then obtain
two complex conjugated branches of eigenvalues
$\lambda_{1,2}(\alpha)=\mu(\alpha)\pm i\omega(\alpha),$ such that
$\lambda_{1,2}(\alpha^*)=\pm i\omega^*.$ If
$\frac{d\mu}{d\alpha}(\alpha^*)\neq 0$ and the first Lyapunov
coefficient of the reduced on the center manifold at $\alpha^*$
equation is different from zero, then a Hopf bifurcation takes place
in the center manifold. The sign of the first Lyapunov coefficient
gives the stability properties of the solution at $\alpha^*$ and of
the periodic solution that occur by Hopf bifurcation. If the first
Lyapunov coefficient is equal to zero, then a degenerated Hopf
bifurcation takes place at $\alpha^*.$

The construction of an approximation of the center manifold and the
computation of the first Lyapunov coefficient (and thus of the
normal form of the reduced equation)  at a Hopf bifurcation point
constitute the object of another paper of ours, \cite{AI-Hopf}.

\textbf{b)} The case $r|\delta+B|=1,$ corresponds to the case
$a_1=1$ from the paper of Hayes, \cite{Ha}. In this case there
always are eigenvalues with either positive or zero real part. The
case of eigenvalues with zero real part (and all other with negative
real part) corresponds to the case $a_2=-1$ of \cite{Ha}. By using
the relations between $a_1,\,a_2$ and $p,\,q,$ (these are
$a_1=-pr,\,a_2=qr$) we find $pr=-1,\,qr=-1,$ hence
$\delta+B=kB=-1/r.$ The first equality implies, as above,
$\frac{\beta_0}{\delta}(k-1)=1$ and it can not be satisfied in the
zone of the parameters that we consider here. It follows that when
$r|\delta+B|=1,$ the solution $x_2$ is unstable.

\textbf{Remark.} Assume that we vary $r$ and keep all the other
parameters fixed. The conditions (\ref{A}) or (\ref{B}) for $r$ are
not as simple as they seem, because $B$ is itself a function of $r$
(being a function of $k$). Let us consider the function
\begin{equation}\label{h}
g(r)=T^{-1}(-(\delta+B(r))r)-\arccos\left(\frac{\delta+B(r)}{k(r)B(r)}\right).
\end{equation}

If for a certain $r^*$ we have $g(r^*)=0$ (that is the condition for
the change of stability), in order to find whether a value $r_1$ in
a neighborhood of $r^*$
 is in the stability zone or not, we have
to know the sign of $g(r_1)$, hence we have to study the monotony
properties of function $g$ in a neighborhood of $r^*$.

\textbf{I.B.} Since here  $B<0,\,B+\delta >0$ we must have

\[ \frac{\beta_0}{\delta}(k-1)>\frac{n}{n-1},\]
\begin{equation}\label{Bplusdpoz}\frac{\beta_0}{\delta}(k-1)(n-k)<n.
\end{equation}
The sufficient condition of local stability is condition (\ref{B}),
with $B$ given by (\ref{B1}).

 A point in the parameter space, satisfying
\[r=\frac{\arccos((\delta+ B)/kB)}{\omega_0}
\]
is  a Hopf bifurcation point, if, with the notations from
\textbf{I.A.a}, $\frac{d\mu}{d\alpha}(\alpha^*)\neq 0$ and the first
Lyapunov coefficient of the reduced on the center manifold at
$\alpha^*$ equation is different from zero. The stability of the
solution in this case is given by the sign of the first Lyapunov
coefficient.

\textbf{Remark.} If $B<0,\,\delta+B=0,$ the solution is stable if
and only if $-kB r<\pi/2$ while, for this case, the point $kB r
=-\pi/2$ is a Hopf bifurcation point.

\textbf{II.} $B>0.$ If $n-1<0,$ then $B>0.$ If $n-1>0,$ then $B>0$
is equivalent to
\begin{equation}\label{condB}
\frac{\beta_0}{\delta}(k-1)<\frac{n}{n-1}.
\end{equation}

In this situation, all the eigenvalues have negative real part if
and only if $kB<\delta+B.$ This inequality is equivalent to
\[\frac{k-1}{\delta}B<1\,\Leftrightarrow\,
\frac{1}{A}[n-(n-1)A]<1\,\Leftrightarrow\,A>1\,\Leftrightarrow\frac{\beta_0(k-1)}{\delta}>1,
\]
and this last inequality is already imposed (by the condition
$x_2>0$. Hence in the case $B>0,$ $x_2$ is stable.

\section{Comments on the stability results in \cite{PM-M}, \cite{PM-B-M}}

In order to compare our results with those of \cite{PM-M},
\cite{PM-B-M}, we define, for $n>1,$
\[r_n :=-\frac{1}{\gamma}\ln \left\{\frac{1}{2}\left(\frac{\delta}{\beta_0}\frac{n}{n-1}+
1\right)\right\}
\]
and remark that $r_n>0\,\Leftrightarrow\,\frac{n}{n-1}\delta<
\beta_0.$ Also for $n>1,$ relation (\ref{condB}) implies that
$B>0\,\Leftrightarrow\,r>r_n$. This last condition is trivially
accomplished if $r_n\leq 0$ which is equivalent to
$\frac{n}{n-1}\delta> \beta_0.$

With these remarks we get the following situations for the sign of
$B$.

I. If $n<1$ then $B>0.$

II. If $n>1$ and $\frac{n}{n-1}\delta> \beta_0$ then $B>0.$

III. If $n>1$ and $\frac{n}{n-1}\delta< \beta_0$ then
\[B>0\,\,\mathrm{for}\,\,r_n<r<r_{max},
\]
\[B<0\,\,\mathrm{for}\,\,0<r<r_n.
\]

This discussion allows us to follow the results of \cite{PM-M},
\cite{PM-B-M} (the delay is there denoted by $\tau$). Those results
have the following weak points.

1. The results in \cite{PM-M} are presented for both equilibrium
points simultaneously, and this leads to imprecisions. As example,
the affirmation at point (1) in \cite{PM-M}, pg. 238, is not true
for $x_1=0.$  Actually, the characteristic equation for this
equilibrium point does not depend on $n$ and thus for this point the
condition $n\in [0,1]$ is irrelevant. The condition of stability for
this point does not depend on $n$. The ambiguity induced by using
the plural "solutions" persists also at point (2) of \cite{PM-M},
pg. 238, leading to misunderstandings since the conclusions there
can not refer to $x_1$, as is seen from our Subsection 2.1.

2. In the case $B<0$, the sign of $\delta+B$  is not considered in
\cite{PM-M}, pg. 238. To express the results of the analysis of the
sign of $\delta+B$ in terms of $r$ is a little more difficult since
inequalities (\ref{Bplusdneg}) and (\ref{Bplusdpoz}) contain second
degree terms in $k$. However the cases $B+\delta>0$ and $B+\delta<0$
are different in conclusions and they can not be eluded.

3. The conclusions in \cite{PM-M}, pg. 238, (2), b) seem to refer to
 the case $B<0,\,\delta+B>0,$ but even for this case the result therein is not correct,
 since there the stability condition is
\[r<\frac{\arccos((B+\delta)/kB)}{\sqrt{(kB)^2-(\delta+B)^2}}\] instead of
\[r<\frac{\arccos((B+\delta)/kB)}{\omega_0}\] with $\omega_0$ defined in
(\ref{om0}), as it is correct (condition (\ref{B})).

Since, in general,  $\omega_0\neq \sqrt{(kB)^2-(\delta+B)^2}$
(equality holds, for $B<0,$  only when relation \eqref{r-H} holds),
it is obvious that the domain of stability found in \cite{PM-M} is
not correct (not even for the case $B<0,\,\delta+B>0$).

In \cite{PM-B-M}, pgs. 316-317, the results are basically the same
as in \cite{PM-M}, excepting the fact that it seems that the
discussion refers only to $x_2$  (but the plural "solutions" is used
again). However, the observations from 2. and 3. above remain valid
for \cite{PM-B-M} also.

\bigskip

Author's address:

"Gh. Mihoc - C. Iacob" Institute of Mathematical Statistics\\
and Applied Mathematics of the Romanian Academy,\\
Calea 13 Septembrie, no. 13, 050711, \\
Bucharest, Romania.

e-mail: anca\_veronica\_ion@yahoo.com

\end{document}